\documentclass[12pt,reqno]{amsart}

\DeclareMathOperator{\chr}{char}

\usepackage{amscd,amsmath,amssymb,latexsym}

\begin{document}

\newtheorem{theorem}{Theorem}
\newtheorem{corollary}[theorem]{Corollary}
\newtheorem{proposition}{Proposition}
\newtheorem{lemma}{Lemma}

\theoremstyle{definition}
\newtheorem*{definition*}{Definition}
\newtheorem*{remark*}{Remark}

\newcommand{\Ec}{\mathcal{E}}
\newcommand{\Fix}{\mathrm{Fix}}
\newcommand{\Fp}{\mathbb{F}_p}
\newcommand{\Gal}{\mathrm{Gal}}
\newcommand{\Jc}{\mathcal{J}}
\newcommand{\Z}{\mathbb{Z}}
\newcommand{\lra}{\longrightarrow}

\newcommand{\mo}{\mathopen\langle}
\newcommand{\mc}{\mathclose\rangle}

\parskip=10pt plus 2pt minus 2pt

\title[Galois Embedding Problems with Cyclic Quotient]
{Galois Embedding Problems with Cyclic Quotient of Order $p$}

\author[J\'{a}n Min\'{a}\v{c}]{J\'an Min\'a\v{c}$^{*\dagger}$}
\address{Department of Mathematics, Middlesex College,
\ University of Western Ontario, London, Ontario \ N6A 5B7 \
CANADA}
\thanks{$^*$Research supported in part by the Natural Sciences and
Engineering Research Council of Canada grant R0370A01, as well as
by the special Dean of Science Fund at the University of Western
Ontario.}
\thanks{$^\dag$Supported by the Mathematical Sciences Research
Institute, Berkeley}
\email{minac@uwo.ca}
\urladdr{http://www.math.uwo.ca/minac.html}

\author[John Swallow]{John Swallow$^\ddag$}
\address{Department of Mathematics, Davidson College, Box 7046,
Davidson, North Carolina \ 28035-7046 \ USA}
\thanks{$^\ddag$Research supported in part by National Security
Agency grant MDA904-02-1-0061.}
\email{joswallow@davidson.edu}
\urladdr{http://www.davidson.edu/math/swallow/}

\begin{abstract}
Let $K/F$ be a cyclic field extension of odd prime degree. We
consider Galois embedding problems involving Galois groups with
common quotient $\Gal(K/F)$ such that corresponding normal
subgroups are indecomposable $\Fp[\Gal(K/F)]$-modules.  For these
embedding problems we prove conditions on solvability, formulas
for explicit construction, and results on automatic realizability.
\end{abstract}

\date{April 5, 2004}

\maketitle

\section*{Introduction}

Let $p>2$ be a prime, and suppose that $K/F$ is a cyclic field
extension of degree $p$. In this paper we consider Galois
embedding problems involving Kummer extensions of $K$ of degree
$p^n$ that are Galois over $F$, and we establish new automatic
realizability results, whereby the solvability of one Galois
embedding problem implies the solvability of another. (See
e.~g.~\cite[Section~5]{GrSmSw} for some automatic realizations of
$2$-groups as Galois groups.) We restrict ourselves to the case
$p>2$ because the case $p=2$ is quite simple and does not lead to
new results.

We focus particularly on the case when $F$ contains a primitive
$p$th root of unity.  In fact, this paper is a continuation of
\cite{MS} wherein, under this hypothesis, we classified
$\Fp[\Gal(K/F)]$-modules $K^\times/K^{\times p}$ using arithmetic
invariants attached to $K/F$, and the investigations there were
motivated by the embedding problems solved in this paper. When $F$
is not of this type, we employ a descent argument in the case
$\chr F\neq p$ and Witt's Theorem in the case $\chr F=p$ to extend
our results to arbitrary fields.

When $F$ contains a primitive $p$th root of unity, we additionally
provide explicit solutions of some Galois embedding problems, and
we show that these formulas are natural and quite transparent
consequences of our method.  For most of these embedding problems,
explicit solutions were not previously known.  For others, such as
the example of Section~\ref{se:exmain}, our methods yield an
explanation of explicit solutions determined previously via ad hoc
methods.

In Section~\ref{se:exmain} we present a motivating example and our
Main Theorem on automatic realizability and explicit solution. In
Section~\ref{se:prelims} we introduce notation and results in
preparation for Section~\ref{se:profunity}, where we give
conditions and explicit solutions for a class of embedding
problems under the hypothesis that a primitive $p$th root of unity
lies in the base field.  In Section~\ref{se:arbitrary} we use a
descent argument and Witt's Theorem to establish equivalent
conditions for embedding problems over all fields, and in
Section~\ref{se:prfmain} we prove our Main Theorem.  Although this
paper uses ideas and results developed in \cite{MS} and in
\cite{Wat}, we decided to make our paper largely self-contained,
and hence we make minimal references to results in \cite{MS} and
\cite{Wat}.

\section{Example and Main Theorem}\label{se:exmain}

A simple example serves as a motivating introduction to Galois
embedding problems of the type we will consider.  Assume that $F$
contains a primitive $p$th root of unity $\xi_p$ and
$K=F(\root{p}\of{a})$ is a cyclic extension of degree $p$, and
consider Heisenberg's group $E$, a noncommutative group of order
$p^3$ and exponent $p$. These conditions determine $E$ up to
isomorphism. The center of $E$ is cyclic of order $p$, and we have
the following short exact sequence:
\begin{equation}\label{eq:shortexact}
    1\to \Z/p\Z \to E\to \Z/p\Z \times \Z/p\Z \to 1.
\end{equation}
Now let $L/K$ be an extension Galois over $F$ such that
$\Gal(L/F)\cong \Z/p\Z \times \Z/p\Z$. Then the exact sequence
naturally gives rise to a Galois embedding problem, asking whether
$L$ embeds in a Galois extension ${\tilde L}/F$ with group $E$ and
such that the surjection in the exact sequence is the surjection
of Galois theory.

The obstruction to the solvability of this embedding problem may
be computed as follows. Assume that $K=F(\root{p}\of{a})$ is contained in $L$.
Fix a primitive root $\xi_p$. By Kummer theory there exist elements
$\sigma,\tau\in\Gal(L/F)$ and $b\in F^\times$ such that
$\root{p}\of{a}^{\sigma-1}=\xi_p=\root{p}\of{b}^{\tau-1}$ and
$\root{p}\of{b}^{\sigma-1}=1=\root{p}\of{a}^{\tau-1}$. Then the lifts
of $\sigma$ and $\tau$ in $E$ generate $E$.
It is well-known that the Galois
embedding problem admits a solution if and only if $b\in
N(K^\times)$, where $N$ denotes the norm map from $K$ to $F$.
(See, for instance, \cite[page~161]{JLY}.)

Moreover, if we suppose that $\omega\in K$ satisfies
$N(\omega)=b$, then it has been observed in
\cite[Cor.~p.~523~\&~Thm.~3(A)]{Ma} (see also
\cite[page~161]{JLY}) that all field extensions ${\tilde L}/F$
solving the Galois embedding problem may be written ${\tilde
L}=L(\root{p}\of{f\alpha})$, where $f\in F^\times$ and
$\alpha=\omega^{p-1}\sigma(\omega)^{p-2}\cdots\sigma^{p-2}(\omega)$.

In our Main Theorem we generalize and motivate both the condition
on solvability and the form of the solution.  The condition
implies that a new automatic realizability result holds for fields
containing $\xi_p$, and we extend the automatic realizability
result to all fields $F$. Further generalizations and explicit
solutions appear in Theorems~\ref{th:ecij}, \ref{th:ecij'}, and
\ref{th:arb}.

Observe that in the example above $L$ and ${\tilde L}$ are Kummer
extensions of $K$ of $p$th-power degree that are Galois over $F$,
and the Galois groups $\Gal(L/K)$ and $\Gal({\tilde L}/K)$ are
naturally acted upon by $\Gal(K/F)$. The appropriate context for
our results turns out to be Kummer extensions $L$ of $K$ such that
$\Gal(L/K)$ is an indecomposable $\Fp[\Gal(K/F)]$-module; as we
show later in Proposition~\ref{pr:compositum}, any Kummer
extension of $K$ of degree $p^n$ that is Galois over $F$
decomposes into a compositum of extensions $L/F$ of this type.

Let $F$ be an arbitrary field, and suppose that $K/F$ is a cyclic
extension with Galois group $G=\Gal(K/F)\cong \Z/p\Z$, with
generator $\sigma$.  Let $A=\oplus_{j=0}^{p-1} \Fp\tau^j$ be a
free $\Fp[G]$-module on the generator $\tau$, where $\sigma$ acts
by multiplication by $\tau$. Let $A_i$ be the $\Fp[G]$-submodule
generated by $(\tau-1)^i$. (See Section~\ref{se:prelims} for
details.)  Finally let $\Ec_i$, $1<i\le p$, denote the following
Galois embedding problem:
\begin{equation*}
    \Ec_i\colon \quad 1\to A_1/A_i \to (A/A_i)\rtimes G\to
    (A/A_1)\rtimes G=\Gal(L/F)\to 1.
\end{equation*}
Observe that $A/A_1\cong \Fp$, a trivial $\Fp[G]$-module; hence
$(A/A_1)\rtimes G\cong \Z/p\Z \times \Z/p\Z$. We also assume that
the projection of $(A/A_1)\rtimes G=\Gal(L/F)$ onto $G$ coincides
with the restriction map $\Gal(L/F)\lra G=\Gal(K/F)$. Assume now
that $\chr F\neq p$. Then $[L:F]$ and $[F(\xi_p):F]$ are coprime.
Therefore $\Gal(L/F)$ is naturally isomorphic to
$\Gal(L(\xi_p)/F(\xi_p))$. After identifying these two Galois
groups we set $F(\xi_p,\root{p}\of{b})$ to be the fixed field
of $1\rtimes G$ in $L(\xi_p)$. (Here $b$ is a suitable element
in $F(\xi_p)^\times$.)

Further observe that in the case $i=2$, $(A/A_2)\rtimes G\cong E$.
Hence $\Ec_2$ is precisely the embedding problem in
equation~\eqref{eq:shortexact} above.

In the following theorem we consider the embedding problems
$\Ec_i$ where $i=2,\dots,p$. We prove:

\begin{theorem}[Main Theorem]\label{th:basicconds} \

    {\bf A)} Let $F$ be an arbitrary field.
         Then the following are equivalent:
    \begin{enumerate}
        \item\label{basicitem1} Some $\Ec_i$ is solvable.
        \item\label{basicitem2} Each $\Ec_i$ is solvable.
    \end{enumerate}
    Consequently, if $(A/A_2)\rtimes G$ occurs as a Galois group
    over $F$, then $(A/A_i)\rtimes G$ occurs as well, for all
    $2\le i\le p.$

    {\bf B)} Now assume that $\chr F\neq p$. Then (\ref{basicitem1}) and
    (\ref{basicitem2}) are also equivalent to
    \begin{enumerate}
        \setcounter{enumi}{2}
        \item\label{basicitem3} $b\in N_{K(\xi_p)/F(\xi_p)}
        (K(\xi_p)^\times)$.
    \end{enumerate}

    {\bf C)} Now assume further that $\xi_p\in F$.  Suppose that
    (\ref{basicitem1})--(\ref{basicitem3}) hold, and let
    $\omega\in K^\times$ satisfy $N(\omega)=b$. Suppose $i>2$. Then a solution
    to $\Ec_i$ is given by
    \begin{equation*}
        {\tilde L}=K(\root{p}\of{f\omega^{(\sigma-1)^{p-i}}},
        \root{p}\of{\omega^{(\sigma-1)^{p-i+1}}}, \dots,
        \root{p}\of{\omega^{(\sigma-1)^{p-2}}}),
    \end{equation*}
    $f\in F^\times$. If $i=2$ then a solution to $\Ec_2$
    is given by $\tilde L=K(\root{p}\of{f\omega^{(\sigma-1)^{(p-2)}}})$.
    Moreover, all solutions of $\Ec_i$ arise in this way.
\end{theorem}

In particular, we have the following automatic realization of
Galois groups: \textit{if $E=(A/A_2)\rtimes G$ is a Galois group
over $F$, $\Fp[G]\rtimes G$ is a Galois group over $F$}.

The explicit construction result in the theorem says that in the
case $\xi_p\in F$, solutions of $\Ec_i$ are parameterized by
$\omega$ with $N(\omega)=b$ and $f\in F^\times$. Note that in
$\Fp[G]$ we have the identity
\begin{equation*}
    (\sigma-1)^{p-2} = (p-1) + (p-2) \sigma + \dots +
    \sigma^{p-2},
\end{equation*}
so the construction of ${\tilde L}$ in the theorem above is
equivalent to that of \cite[page~161]{JLY} in the case $i=2$.

\section{Preliminaries}\label{se:prelims}

In this Section and Section~\ref{se:profunity} we assume that $F$
is a field containing a primitive $p$th root of unity $\xi_p$,
$K=F(\root{p}\of{a})$, and $G=\Gal(K/F)\cong \Z/p\Z$.  We let
$\sigma$ denote the generator of $G$ such that $\root{p}
\of{a}^{\sigma-1}= \xi_p$. Since $\sigma-1$ is used frequently, we
use the abbreviation $\rho:=\sigma-1$. All modules and Galois
extensions will be acted upon on the left by their respective
groups, even though we will use exponential notation to denote
Galois action on fields.  We denote by $F^\times$ the
multiplicative group of a field $F$, and we write $N=N_{K/F}$ for
the norm map from $K$ to $F$. For a subset $S$ of an $\Fp$-module
$V$ we denote by $\mo S \mc$ the $\Fp$-span of $S$ in $V$.

\subsection{$\Fp[G]$-modules}\

Let $A=\oplus_{j=0}^{p-1} \Fp\tau^j$ be a free $\Fp[G]$-module
on the generator $\tau$, where $\sigma$ acts by multiplication by
$\tau$.  There are $p$ quotient modules $A/A_i,i=1,\dots,p$ of $A$
where for $i<p$,
\begin{equation*}
    A_i=\mo (\tau-1)^i, (\tau-1)^{i+1}, \dots, (\tau-1)^{p-1} \mc,
    \quad \text{and} \quad A_p=\{0\}.
\end{equation*}
These quotients are all cyclic and together form a complete set of
indecomposable $\Fp[G]$-modules.  Each $A/A_i$ is of dimension $i$
as a vector space over $\Fp$.  We call this dimension the
\emph{length}, and denote the length of a cyclic $\Fp[G]$-module
$M$ by $l(M)$, because we have the following
criterion for $l(M)$, where
$M$ is a cyclic $\Fp[G]$-module generated by $m$: $l(M)=i$ such
that $\rho^im=0$, $\rho^{i-1}m\neq 0$.  Moreover,
such a cyclic module $M$ of length $l$ contains precisely one
submodule of each length $1\le j\le l$: $M_i=\mo
\rho^{l-j}m,\dots,\rho^{l-1}m\mc$.

For each $i\in\{1,\dots,p\}$ we pick a basis
$\{1,\overline{\tau-1},\dots,(\overline{\tau-1})^{i-1}\}$ of
$A/A_i$ consisting of images of $1,\tau-1,\dots,(\tau-1)^{i-1}$.
We define an $\Fp$-linear map $\lambda\colon A/A_i\lra\Fp$ by
$\lambda(f_0+f_1(\overline{\tau-1})+\dots+
f_{i-1}(\overline{\tau-1})^{i-1})=f_{i-1}$, where
$f_k\in\Fp,k=0,\dots,i-1$. Observe that $\ker(\lambda)$ contains
no nonzero ideal of $A/A_i$. Then $B(a,b):=\lambda(ab)$ for each
$a,b\in A/A_i$ defines a nonsingular, symmetric bilinear form
$B\colon A/A_i\times A/A_i\lra\Fp$. Thus $A/A_i$ is a symmetric
algebra. (See \cite[page~442]{Lam}.) Further we have
$B(a^{\sigma},b^{\sigma^{-1}})=B(a,b)$ for each $a,b\in A/A_i$ and
our bilinear form $B$ induces a $G$-equivariant isomorphism
between $A/A_i$ and its dual.

\subsection{Groups}\

In this section we classify the groups of interest in this paper
and the surjections among them. For $e\in \Fp$, let $B_{i,e}$ be
the group extension of $A/A_i$ by $G$ with $\tilde\sigma^p=
e(\overline{\tau-1})^{i-1}$. Here $\tilde\sigma$ is a lift in
$B_{i,e}$ of $\sigma\in G$. Note that for $e=0$, $B_{i,0}=(A/A_i)
\rtimes G$. First we consider the equivalence classes of these
groups.

\begin{lemma}\label{le:group} (See \cite[Theorem 2]{Wat})\
    \begin{enumerate}
        \item If $H$ is a group with a normal subgroup isomorphic
        to $A/A_i$ as a $G$-module, with quotient group $G$, then
        $H=B_{i,e}$ for some $e$.
        \item For fixed $1\le i<p$, all $B_{i,e}$, $e\neq 0$, are
        isomorphic, and these groups are not isomorphic to
        $B_{i,0}=(A/A_i) \rtimes G$.
        \item For $i=p$, all $B_{i,e}$ are isomorphic to
        $B_{p,0}\cong \Fp[G]\rtimes G$.
    \end{enumerate}
\end{lemma}

The Galois embedding problems in this paper consist of embedding
an extension $L/F$ with group $B_{j,e'}$ in an extension with
strictly larger group $B_{i,e}$. We are interested in all
surjections $B_{i,e}\to B_{j,e'}$ for which the kernel lies in
$A/A_{i}\subset B_{i,e}$ and which are induced by the projection
of $B_{i,e}$ on its quotient. We call these $G$-surjections.

\begin{lemma}\label{le:surjections}
    The $G$-surjections in the set of groups $\{B_{i,e}\}_{i\ge
    1}$ are precisely
    \begin{equation*}
        B_{i,e}\to B_{j,0}, \quad i>j\ge 1, \ e\in \Fp, \quad
        \text{with kernel\ } A_j/A_i.
    \end{equation*}
\end{lemma}

\begin{proof}
    Considering the dimensions of $A/A_i$ and $A/A_j$, if
    $B_{i,e}\to B_{j,e'}$ is a $G$-surjection then $i>j$.  Now a
    surjection of $G$-modules $A/A_i$ to $A/A_j$ must have as
    kernel an $\Fp[G]$-submodule of $A/A_i$ of $\Fp$-rank $i-j$.
    But since $A/A_i$ is cyclic, there is precisely one such
    submodule, namely $A_j/A_i$.  Hence $(\overline{\tau-1})^k$
    lies in the kernel for all $j\leq k<i$.  In particular, the
    kernel must contain $e(\overline{\tau-1})^{i-1}$, which is
    $\tilde\sigma^p$ in $B_{i,e}$. Therefore $\tilde\sigma\in
    B_{i,e}$ is sent to some lift $\hat\sigma\in B_{j,e'}$ of
    $\sigma\in G$ and hence $\hat\sigma^p=1$ in $B_{j,e'}$, or
    $e'=0$.
\end{proof}

We list some characteristics of the groups $B_{i,e}$. Each
$B_{i,e}$ has order $p^{i+1}$, nilpotent index $i$, and rank (the
smallest number of generators) $2$. The exponent of $B_{i,0}$ is
$p$, and the exponent of $B_{i,e}, e\neq 0$ is $p^2$.  The
Frattini subgroup $\Phi(B_{i,e})$ of $B_{i,e}$ is $A_1/A_i \cong
(\Z/p\Z)^{i-1}$.  Finally, we have presentations
\begin{multline*}
    B_{i,0} = \langle  \sigma, \ \ \{\tau_j\}_{j=0}^{i-1}
    \quad : \\ \sigma^p=\tau_j^p=[\sigma,\tau_{i-1}]=1; \ \
    \text{for $j<i-1$}, \ [\sigma,\tau_j]=\tau_{j+1} \rangle
\end{multline*}
and, for $e\not\equiv 0\bmod p$,
\begin{multline*}
    B_{i,e} = \langle \sigma, \ \ \{\tau_j\}_{j=0}^{i-1}
    \quad : \quad \sigma^p=\tau_{i-1}^e; \\
    \tau_j^p=[\sigma,\tau_{i-1}]=1; \ \
    \text{for $j<i-1$},\
    [\sigma,\tau_j]=\tau_{j+1} \rangle.
\end{multline*}

\subsection{Extensions and Submodules, $\xi_p\in F$}\

Now let $J$ denote the $\Fp[G]$-module $J:=K^\times/K^{\times p}$.
We denote elements of $J$ by $[\gamma]$, $\gamma\in K^\times$. Let
$J_{i}$ be the kernel of the endomorphism $(\sigma-1)^i$ and let
$M_{\gamma}$ be the cyclic submodule of $J$ generated by
$[\gamma]$. Then $[\gamma]\in J_i$ if and only if $l(M_\gamma)\le
i$.

We denote by $M\leftrightarrow L_M$ the Kummer correspondence over
$K$ of subspaces $M$ of the $\Fp$-vector space $J$ and abelian
exponent $p$ extensions $L_M$ of $K$:
\begin{equation*}
    M = (L_M^{\times p}\cap K^{\times})/{K^{\times p}}
    \quad \leftrightarrow \quad
    L_M = K(\root{p}\of{\gamma} \colon [\gamma]\in M).
\end{equation*}
Set $C=\Gal(L_M/K)$. Then $M$ and $C$ are dual $G$-modules and the
canonical duality $\mo m,c \mc:=c(\root{p}\of{m})/\root{p}\of{m}$
of $M$ and $C$ is $G$-equivariant. (See \cite[pages~134 and
135]{Wat}.) The following proposition rephrases the results in
\cite[page~135]{Wat} in our notation.

\begin{proposition}\label{pr:Kummer}
    Under the Kummer correspondence above,
    \begin{enumerate}
        \item\label{item1} $L_M$ is Galois over $F$ if and only if
        $M$ is an $\Fp[G]$-submodule of $J$.
        \item\label{item2} The following are equivalent:
        \begin{enumerate}
            \item\label{it2a} $L_M$ is the Galois closure, over
            $F$, of $K(\root{p}\of{\gamma})$ for some $\gamma\in
            K^\times$;
            \item\label{it2b} $M=M_{\gamma}$ for some $\gamma\in
            K^\times$;
            \item\label{it2c} $\Gal(L_M/K)\cong A/A_i$, as
            $G$-modules, for some $i$;
            \item\label{it2d} $\Gal(L_M/F)\cong
            B_{i,e}$, as $G$-extensions, for some $i$ and $e$.
        \end{enumerate}
        If these conditions hold, then $i=l(M)$ and
        \begin{equation*}
            L_M = K(\root{p}\of{\gamma},
            \root{p}\of{\gamma^{\rho}}, \dots,
            \root{p}\of{\gamma^{\rho^{i-1}}}).
        \end{equation*}
    \end{enumerate}
\end{proposition}

\begin{proof}
    Because $L_M$ is Galois if and only if each automorphism of
    $K$ extends to an automorphism of $L_M$, item (\ref{item1})
    and (\ref{it2a})$\Leftrightarrow$(\ref{it2b}) follow. That
    (\ref{it2c})$\Leftrightarrow$(\ref{it2d}) follows from
    Lemma~\ref{le:group}.

    Suppose (\ref{it2b}) holds. Then $M\cong A/A_i$ for some
    $i\in\{1,\dots,p\}$ and $\Gal(L_M/K)$ is a $G$-equivariant
    dual of $M$. Since $M$ is a $G$-equivariant self-dual module,
    we see that $\Gal(L_M/K)$ and $A/A_i$ are $G$-isomorphic and
    (\ref{it2c}) follows.

    Suppose now that (\ref{it2c}) holds. Then again using the
    $G$-equivariant self-duality of $A/A_i$ and Kummer theory, we
    see that $M$ must be a cyclic module $M_\gamma$ for some
    $\gamma\in K^{\times}$. Hence (\ref{it2b}) follows.

    The presentation of $L_M$ follows from the fact that a cyclic
    $\Fp[G]$-module $M$ generated by $m$ is generated over $\Fp$
    by $\{\rho^k(m)\}_{k=0}^{l(M)-1}$.
\end{proof}

We can now prove

\begin{proposition}\label{pr:compositum}
    Let $L/K$ be a finite Kummer extension of $p$th-power degree
    which is Galois over $F$.  Then $L$ is a compositum of
    finitely many Galois closures, over $F$, of extensions of the
    form $L_\gamma=K(\root{p}\of{\gamma}),$ $\gamma\in K^\times$.
\end{proposition}

\begin{proof}
    The extension $L/K$ corresponds to an $\Fp[G]$-submodule $M$
    of $J$. Since $M$ is finite, it is decomposable into a direct
    sum of finitely many indecomposable $\Fp[G]$-modules $M_j$.
    Each indecomposable $\Fp[G]$-module $M_j$ is isomorphic to
    some $A/A_i$ and is hence cyclic.  By
    Proposition~\ref{pr:Kummer} (\ref{item2}), these submodules
    correspond to Galois closures over $F$ of extensions
    $L_{\gamma}=K(\root{p}\of{\gamma})$.  The submodule of $J$
    generated by each of the indecomposables $M_j$ then
    corresponds to the compositum of the $L_\gamma$, and we are
    done.
\end{proof}

\subsection{The index}\label{se:index}\

The following homomorphism appears in a somewhat different form in
\cite[Theorem 3]{Wat}:

\begin{definition*}
    The \emph{index} $e([\gamma])\in \Fp$ for
    $[\gamma]\in J_{p-1}$ is defined by
    \begin{equation*}
        \xi_p^{e([\gamma])} =
        \left(\root{p}\of{N_{K/F}(\gamma)}\right)^{\rho}.
    \end{equation*}
\end{definition*}
The index is well-defined, as follows.  First, since
\begin{equation}\label{eq:equiv}
    1+\sigma+\dots+\sigma^{p-1}=(\sigma-1)^{p-1}=\rho^{p-1}
\end{equation}
in $\Fp[G]$, $[N(\gamma)]=[\gamma]^{\rho^{p-1}}$, which is the
trivial class $[1]$ by the assumption $[\gamma]\in J_{p-1}$, and
as a result $\root{p}\of{N(\gamma)}$ lies in $K$ and is acted upon
by $\sigma$.  Observe further that $e([\gamma])$ depends neither
on the representative $\gamma$ of $[\gamma$] nor on the particular
$p$th root of $N(\gamma)$. Also the index function $e$ above is a
group homomorphism from $J_{p-1}$ to $\Fp$. Therefore the
restriction of $e$ to any $M_{\gamma}$ is either trivial or
surjective.

We show that the index is trivial for any $[\gamma]$ in the image
of $\rho$:
\begin{equation*}
    \xi_p^{e([\gamma]^{\rho})} = \root{p}\of{
    N(\gamma^{\sigma-1})}^{\rho} = \root{p}\of{1}^{\rho} = 1,
\end{equation*}
or $e([\gamma]^{\rho})=0$.

\begin{lemma}\label{le:groupoflm} (See \cite[Theorem 2]{Wat}.)
    Let $[\gamma]\in J$ and $M=M_\gamma$.
    \begin{enumerate}
        \item If $l(M)<p$ and $e=e ([\gamma])$ then
        $\Gal(L_M/F)\cong B_{i,e}$.
        \item If $l(M)=p$ then $\Gal(L_M/F)\cong B_{p,0}$.
    \end{enumerate}
\end{lemma}

\begin{proof}
    The second item follows from Proposition~\ref{pr:Kummer} and
    Lemma~\ref{le:group}.  The fact that $\Gal(L_M/F)\cong
    B_{i,e}$ for some $e\in \Fp$ follows in the same
    manner.  Therefore it remains only to show that
    $\Gal(L_M/F)\cong B_{i,e([\gamma])}$.

    Let $\tilde\sigma$ denote a pullback of $\sigma\in G$ to
    $\Gal(L_M/F)$. Then $\tilde\sigma^{p}$ lies in $Z(\Gal(L_M/F))
    \cap \Gal(L_M/K)$. (Here $Z(\Gal(L_M/F))$ means the center of
    $\Gal(L_M/F).)$ Recall that using Kummer theory and the
    $G$-equivariant self-duality of $A/A_i$ we may identify
    $\Gal(L_M/K)$ with $A/A_i$. Adopting this identification we
    pick a basis $\{1,\overline{\tau-1},\dots,
    (\overline{\tau-1})^{i-1}\}$ of the $\Gal(L_M/F)$ dual with
    $\{[\gamma]^{(\sigma-1)^{i-1}}, \dots,[\gamma]^{\sigma-1},
    [\gamma]\}$ with respect to Kummer pairing. Under our
    identification, $\tilde\sigma^{p}$ lies in the $G$-invariant
    submodule of $A/A_{i}$, which is $\mo(\overline{
    \tau-1})^{i-1}\mc$. Observe that $(\overline{\tau-1})^{i-1}$
    sends $\root{p}\of{\gamma}$ to $\xi_p\root{p}\of{\gamma}$. If
    $\tilde\sigma^{p}=e(\overline{\tau-1})^{i-1}$ then
    \begin{equation*}
        (\root{p}\of{\gamma})^{(\tilde\sigma^{p}-1)} = \xi_{p}^{e}.
    \end{equation*}
    Therefore
    \begin{equation*}
        \root{p}\of{\gamma}^{(\tilde\sigma^{p}-1)} =
        \root{p}\of{\gamma}^{(1+\tilde\sigma+
        \dots+\tilde\sigma^{p-1}) (\tilde\sigma-1)} =
        \left(\root{p}\of{
        N_{K/F}(\gamma)}\right)^{(\tilde\sigma-1)}=
        \xi_p^{e([\gamma])}.
    \end{equation*}
\end{proof}

We characterize elements of $J$ fixed by $\sigma$ and of trivial
index with the following

\begin{lemma}\label{le:j1triv}
    If $[\gamma]\in J_1$ and $e([\gamma])=0$ then there exists
    $f\in F^\times$ such that $[\gamma]=[f]$.
\end{lemma}

\begin{proof}
    By \cite[Remark 2]{MS}, we have the following short exact
    sequence:
    \begin{equation*}
        0\to \mo [a] \mc \xrightarrow{i} F^\times/F^{\times p}
        \xrightarrow{\epsilon} J_1 \xrightarrow{N} \mo [a]\mc,
    \end{equation*}
    where $\mo [a]\mc$ is the subgroup of $F^\times/F^{\times p}$
    generated by $[a]\in F^\times/F^{\times p}$, $i$ is the
    inclusion map, $\epsilon$ is the natural homomorphism induced
    by the inclusion map $F^\times \to K^\times$, and $N$ is the
    map induced by the norm map from $K$ to $F$. Now
    $e([\gamma])=0$ implies that $[\gamma]$ is in the kernel of
    the surjection $N$ above, and we are done.
\end{proof}

We will also need a lemma on the smallest lengths of cyclic
submodules of $J$ generated by an element $[\gamma]$ with
nontrivial index.  Let $\Upsilon=1$ if $\xi_p\in N(K^\times)$ and
$\Upsilon=0$ otherwise. In the proof of the next lemma we refer to
\cite{A} only for the sake of convenience. One can use basic
Kummer theory instead.

\begin{lemma}\label{le:smlgthnontriv}\
    \begin{enumerate}
        \item If $\Upsilon=1$ then there exists $\delta\in
        K^\times$ such that $[\delta]\in J_1$ and $e([\delta])
        \neq 0$. These are precisely the $\delta$ such that
        $K(\root{p}\of{\delta})/F$ is a cyclic extension of
        degree $p^2$.
        \item If $\Upsilon=0$ then $[\root{p}\of{a}]\in
        J_2\setminus J_1$, $e([\root{p}\of{a}])\neq 0$, and
        $e([\gamma])=0$ for all $[\gamma]\in J_1$.
    \end{enumerate}
\end{lemma}

\begin{proof}
    By \cite[Theorem 3]{A}, $\Upsilon=1$ if and only if
    $K/F$ embeds in an extension $L=K(\root{p}\of{\delta})$ Galois
    over $F$ with group $\Z/p^2\Z \cong B_{1,e}$, $e\neq 0$.  By
    Proposition~\ref{pr:Kummer} and Lemma~\ref{le:groupoflm},
    then, $\Upsilon=1$ if and only if there exists $[\delta]\in
    J_1$ with $e([\delta]) \neq 0$.  This proves the first
    statement.

    Assume now that $\Upsilon=0$. We have $[\root{p}\of{a}]\in
    J_2$, since $[\root{p}\of{a}]^{\rho} = [\xi_p] \in J_1$, and
    we calculate $e([\root{p}\of{a}])=1$. Since $\Upsilon=0$,
    $[\xi_p]\neq [1]$ in $J_1$ and therefore $[\root{p}\of{a}]
    \not\in J_1$. Now consider any $[\gamma]\in J_1$. Then
    $L=K(\root{p}\of{\gamma})$ is Galois over $F$ and since
    $\Upsilon=0$ we see from \cite[Theorem 3]{A} that $\Gal(L/F)$
    is $B_{1,0}\cong\Z/p\Z\times\Z/p\Z$. Hence $e([\gamma])=0$ and
    the second statement is proved.
\end{proof}

Finally, we introduce a variant of \cite[Lemma 1]{MS} for
submodules generated by elements with trivial index. This is our
key lemma:

\begin{lemma}\label{le:extend}
    Let $[\gamma]\in J$.  Suppose that
    $2\le l(M_\gamma)<p$ and $e([\gamma])=0$.

    Then there exists $[\gamma']\in J$ such that
    \begin{enumerate}
        \item $l(M_{\gamma'})=l(M_\gamma)+1$.
        \item $[\gamma']^{\rho^2}=[\gamma]^{\rho}$.
        \item The fixed elements $M_{\gamma}^G$ of $M_{\gamma}$
        under $G$ coincide with $M_{\gamma'}^G$.
        \item If $l(M_\gamma)<p-1$ then $e([\gamma'])$ is defined
        and has a value of $0$.
    \end{enumerate}
\end{lemma}

\begin{proof}
    Let $c=N\gamma$.  Since $l(M_\gamma)<p$, we have
    $[c]=[\gamma]^{ \rho^{p-1}} = [1]$. Hence $c\in F^\times \cap
    K^{\times p}$.  In fact, $c= a^sf^p$ for some $f\in F$ and
    $s\in \Z$, as follows. Since $c\in K^{\times p}$,
    $F(\root{p}\of{c})\subset K$.  The Kummer extension
    $F(\root{p}\of{c})$ is either $F$ or $K$; if the former, then
    $c\in F^{\times p}$, while if the latter, then by Kummer
    theory $c$ also has the desired form.

    Thus $N\gamma=a^sf^p$ for some $s$ and $f$.  But $e(\gamma)=0$,
    so $p$ divides $s$ and we see that $N\gamma=f^p$ for
    some $f\in F^\times$.  Since
    $N(\gamma/f)=1$, by Hilbert's Theorem 90 there exists a
    $\omega\in K^\times$ such that $\omega^{\sigma-1}=\gamma/f$.
    Then $l(M_\omega)=l(M_\gamma)+1$.

    If $l(M_\gamma)<p-1$ then let $t=e([\omega])$ and set
    $\gamma'=\omega/(a^{t/p})$; otherwise let $t=0$ and set
    $\gamma'=\omega$.

    We compute $[\gamma']^{\rho}= [\xi_p^{-t} \gamma/f]$ and,
    since $\xi_p$, $f\in F^\times$, $[\gamma']^{\rho^2} =
    [\gamma]^{\rho}$, which is nontrivial since $2\le
    l(M_\gamma)$. Hence (1) and (2) follow.

    Now if $l(M_\gamma)<p-1$ then $l(M_{\gamma'})=l(M_\gamma)+1<p$
    and so $\gamma'\in J_{p-1}$.  Then
    $e([\gamma'])=e([\omega])-t=0$. Therefore (4) is valid.

    Finally observe that $M_\gamma^G$ is generated by
    $[\gamma]^{\rho^{l(M_\gamma)-1}}$ as well as
    $[\gamma']^{\rho^{l(M_{\gamma'})-1}}$, which in turn generates
    $M_{\gamma'}^G$. Hence $M_\gamma^G=M_{\gamma'}^G$ and
    therefore (3) follows from (2).
\end{proof}

\section{Embedding Problem Conditions and
Solutions, $\xi_p\in F$}\label{se:profunity}

We consider all embedding problems involving groups $B_{i,e}$,
based on the $G$-surjections determined in
Lemma~\ref{le:surjections}, defining the following embedding
problems for $i>j\ge 1$:
\begin{equation*}
    \Ec_{i,j}(L)\colon \quad 1\to A_j/A_i \to B_{i,0} \to
    (A/A_j)\rtimes G=\Gal(L/F)\to 1.
\end{equation*}
and, for any $e\neq 0$,
\begin{equation*}
    \Ec_{i,j}'(L)\colon \quad 1\to A_j/A_i \to B_{i,e}\to
    (A/A_j)\rtimes G=\Gal(L/F)\to 1.
\end{equation*}
In each case we ask if there exists a Galois extension ${\tilde
L}/F$ containing $L$ such that $\Gal(\tilde{L}/F)\cong B_{i,0}$
or $\Gal(\tilde{L}/F)\cong B_{i,e}$,
and under the identification of $\Gal(\tilde{L}/F)$ with $B_{i,0}$
(or $B_{i,e}$), the surjection $\Gal({\tilde
L}/F)\to \Gal(L/F)$ is identical to the surjection above.

Since $B_{p,0}\cong B_{p,e}$ as $G$-extensions for all $e$, the
embedding problems $\Ec_{p,i}(L)$ and $\Ec_{p,i}'(L)$ are
identical. Moreover note $\Ec_{i,1}(L)=\Ec_i(L)$.
(See page~2 for the discussion of $\Ec_i$.)

For each of these problems, by Proposition~\ref{pr:Kummer}, $L$ is
the Galois closure of $K(\root{p}\of{\gamma})$ for some $\gamma\in
K^\times$.  Hence under the Kummer correspondence $M_{\gamma}
\leftrightarrow L$, and by Proposition~\ref{pr:Kummer},
$l(M_\gamma)=j$.

\begin{theorem}\label{th:ecij}
    Suppose that $\xi_p\in F$.  Let $p\geq i>j\geq 1$ and $L$ be
    the Galois closure of $K(\root{p}\of{\gamma})$ over $F$.

    Then $\Ec_{i,j}(L)$ is solvable if and only if $[\gamma] =
    [\omega]^{\rho^{p-j}}$ for some $\omega\in K^\times$.

    If so, then a solution ${\tilde L}$ to $\Ec_{i,j}(L)$, where
    $i>j+1$, is given by
    \begin{equation*}
        {\tilde L}=L(\root{p}\of{f\omega^{\rho^{p-i}}},
        \root{p}\of{\omega^{\rho^{p-i+1}}}, \dots,
        \root{p}\of{\omega^{\rho^{p-j-1}}}),
    \end{equation*}
    $f\in F^\times$.

    In the case when $i=j+1$ a solution $\tilde{L}$ to
    $\Ec_{j+1,j}$ is given by
    $\tilde{L}=L(\root{p}\of{f\omega^{\rho^{p-i}}})$, $f\in
    F^\times$.

    Moreover, all solutions to $\Ec_{i,j}(L)$ arise in this way if
    one allows $\omega$ to vary over all elements of $K^\times$
    with $[\omega]^{\rho^{p-j}}=[\gamma]$.
\end{theorem}

\begin{proof}
    By Proposition~\ref{pr:Kummer}, there exists a field ${\tilde
    L}$ with $\Gal({\tilde L}/F)\cong B_{i,e}$ for some $i$ and
    $e$ if and only if there exists a cyclic submodule $M_{\beta}$
    of $J$ of length $i$, and in this case we have $M_{\beta}
    \leftrightarrow {\tilde L}$ under the Kummer correspondence.

    Furthermore, by Lemma~\ref{le:groupoflm}, if $i<p$ then
    $\Gal({\tilde L}/F)\cong B_{i,e}$, where $e=e([\beta])$, and
    if $i=p$ then $\Gal({\tilde L}/F) \cong B_{p,0}$. Hence if
    $i<p$ then $\Ec_{i,j}(L)$ is solvable if and only if there
    exists $\beta\in K^\times$ with $e([\beta])=0$,
    $l(M_\beta)=i$, and $M_\beta\supset M_\gamma$. If $i=p$ then
    $\Ec_{p,j}$ is solvable if and only if there exists $\beta\in
    K^\times$ with $l(M_\beta)=p$ and $M_\beta\supset M_\gamma$.

    Now suppose that $[\gamma] = [\omega]^{\rho^{p-j}}$ for some
    $\omega\in K^\times$.  Then let $\beta=\omega^{\rho^{p-i}}$.
    Since $l(M_\gamma)=j$ and $[\gamma]=[\beta]^{\rho^{i-j}}$,
    $l(M_\beta)=i$.  Now if $i=p$ then the condition of the
    previous paragraph is satisfied.  If $i<p$ then $\beta$ is in
    the image of the endomorphism $\rho$, therefore $e([\beta])=0$
    and the condition of the previous paragraph is satisfied.

    Going the other way, suppose that there exists $\beta\in
    K^\times$ with $l(M_\beta)=i$, $M_\beta\supset M_\gamma$, and,
    if $i<p$, $e([\beta])=0$.  Since $M_\gamma$ is the unique
    $\Fp[G]$-submodule of $M_\beta$ of length $j$, $M_\gamma =
    M_{\beta^{\rho^{i-j}}}$. Further since the linear map
    $M_{\beta}\to M_{\gamma}$ defined by
    $[\alpha]\mapsto[\alpha]^{\rho^{i-j}}$ is surjective, there is
    a $[\beta']\in M_{\beta}$ such that we have
    $[\beta']^{\rho^{i-j}}=[\gamma]$.  Moreover,
    $l(M_{\beta'})=l(M_\beta)$ so $M_{\beta'}=M_{\beta}$. In the
    case of $i<p$, because $e$ is trivial on $[\beta]$, then $e$
    is trivial on $M_\beta$ and hence $e([\beta'])=0$.

    If $i=p$ then let $\omega=\beta'$.   Otherwise, by repeated
    application of Lemma~\ref{le:extend}, we may find an
    $\omega\in K^\times$ such that $[\omega]^{\rho^{p-i+1}}
    = [\beta']^{\rho}$.  Then $[\omega]^{\rho^{p-j}}
    =[\beta']^{\rho^{i-j}}= [\gamma]$.

    We now treat the explicit construction of the solution fields.
    Let $M_\beta$ be an $\Fp[G]$-module corresponding to a
    solution field to the embedding problem $\Ec_{i,j}(L)$. Let
    $\beta'$ and $\omega$ be defined as above. Note that
    $[\omega]^{\rho^{p-i+y}}=[\beta']^{\rho^y}$ for all $y\ge 1$.
    Hence $\delta=\omega^{\rho^{p-i}}/\beta'$ satisfies
    $[\delta]^{\rho}=[1]$.  Now $e([\delta])=0$, so by
    Lemma~\ref{le:j1triv}, $[\delta]=[f]$ for some $f\in
    F^\times$. If we have a solution ${\tilde L}$ to
    $\Ec_{i,j}(L)$ with $M_{\beta} \leftrightarrow {\tilde L}$,
    then ${\tilde L}=L(\root{p}\of{\theta}: [\theta]\in M_\beta)$,
    or
    \begin{equation*}
        {\tilde L} = L(\root{p}\of{\beta},
        \root{p}\of{\beta^{\rho}}, \dots,
        \root{p}\of{\beta^{\rho^{i-j-1}}}),
    \end{equation*}
    by Proposition~\ref{pr:Kummer}.  Since $M_{\beta'}=M_\beta$,
    $[\omega]^{\rho^{p-i+y}}=[\beta']^{\rho^y}$ for all $y\ge 1$,
    and $[\beta']=[\omega^{\rho^{p-i}}/f]$, we have
    \begin{equation*}
        {\tilde L}= L(\root{p}\of{f^{-1}\omega^{\rho^{p-i}}},
        \root{p}\of{\omega^{\rho^{p-i+1}}},
        \dots,\root{p}\of{\omega^{\rho^{p-j-1}}}),
    \end{equation*}
    in the case when $i>j+1$ and
    $\tilde{L}=L(\root{p}\of{f^{-1}\omega^{\rho^{p-i}}})$ if
    $i=j+1$ again by Proposition~\ref{pr:Kummer}.

    Finally observe that if we have a solution $\tilde{L}$ to
    $\Ec_{i,j}(L)$ with $M_\beta\leftrightarrow\tilde{L}$, then
    for each $f\in F^\times$ a module $M_{f\beta}$ also
    corresponds to a solution of $\Ec_{i,j}(L)$. Hence in our
    explicit formula for a solution field $\tilde{L}$, any $f\in
    F^\times$ is eligible.
\end{proof}

\begin{theorem}\label{th:ecij'}
    Suppose $\xi_p\in F$.  Let $p\geq i>j\ge 1$ and $L$ be the
    Galois closure of $K(\root{p}\of{\gamma})$ over $F$.

    \begin{enumerate}
        \item\label{it:ecij'1} $\Ec_{i,j}'(L)$, $i>j+1-\Upsilon$
        or $j=p-1$, is solvable if and only if $[\gamma] =
        [\omega]^ {\rho^{p-j}}$ for some $\omega\in
        K^\times$.  If so, then a solution ${\tilde L}$ to
        $\Ec_{i,j}'(L)$ is given by
        \begin{equation*}
            {\tilde L} =
            L(\root{p}\of{f\alpha\omega'^{\rho^{p-i}}},
            \root{p}\of{\omega'^{\rho^{p-i+1}}}, \dots,
            \root{p}\of{\omega'^{\rho^{p-j-1}}}),
        \end{equation*}
        $f\in F^\times$, where in the case $\Upsilon=1$, $\alpha$
        is any element in $K^\times$ with $K(\root{p}\of
        {\alpha})/F$ cyclic of degree $p^2$, and in the case
        $\Upsilon=0$, $\alpha$ is $\root{p}\of{a}$. Furthermore
        $\omega' = \omega^{c_0+c_1\rho+ \dots+ c_{i-1}\rho^{i-1}}$
        for suitable $c_k\in\mathbb{Z}$.
        \bigskip
        \item\label{it:ecij'2} $\Ec_{j+1,j}'(L)$, $\Upsilon=0$, is
        solvable if and only if
        $[\gamma]=[\xi_p]^e[\omega]^{\rho^{p-j}}$ for some
        $\omega\in K^\times$ and $e\not\equiv 0\bmod p$. If so,
        then a solution ${\tilde L}$ to $\Ec_{i,j}'(L)$ is given
        by
        \begin{equation*}
            {\tilde L}= L(\root{p}\of{fa^{e/p}
            \omega^{\rho^{p-j-1}}}), \ \ f\in F^\times.
        \end{equation*}
    \end{enumerate}
    Moreover, all solutions to $\Ec_{i,j}'$ arise in the way
    described above.
\end{theorem}

Note that the two parts of the theorem overlap when $i=p$,
$j=p-1$, and $\Upsilon=0$.

\begin{proof}
    We begin in the same manner as the previous proof: if $i<p$
    then $\Ec_{i,j}'(L)$ is solvable if and only if there exists
    $\beta\in K^\times$ with $e([\beta])\neq 0$, $l(M_\beta)=i$,
    and $M_\beta\supset M_\gamma$.  If $i=p$ then $\Ec_{p,j}'$ is
    solvable if and only if there exists $\beta\in K^\times$ with
    $l(M_\beta)=p$ and $M_\beta\supset M_\gamma$.

    We first treat the conditions on $[\gamma]$ that are equivalent
    to solvability.  Then we consider the explicit presentations
    of the solution fields.

    In the case $i=p$, since $\Ec'_{p,j}=\Ec_{p,j}$, the condition
    on $[\gamma]$ is the same as the condition on $[\gamma]$ for
    the solvability of $\Ec_{p,j}$ determined in the previous
    theorem. This gives us the condition in part~\ref{it:ecij'1}.
    Now if additionally $j=p-1$, consider the condition of
    part~\ref{it:ecij'2}: $[\gamma]= [\xi_p]^e[\omega']^
    {\rho^{p-j}}=[\xi_p]^e[\omega']^{\rho}$ for $e\not\equiv 0
    \bmod p$. If this condition holds, $\omega=a^{-e/p}\omega'$
    satisfies the condition $[\gamma]=[\omega]^{
    \rho^{p-j}}=[\omega]^{\rho}$ of part~\ref{it:ecij'1}.
    Conversely, if the condition of part~\ref{it:ecij'1} holds,
    set $\omega'=a^{1/p}\omega$ and observe that the condition of
    part~\ref{it:ecij'2} holds with $e=1$.

    Now suppose $i<p$ and $\Ec_{i,j}'(L)$ is solvable with field
    ${\tilde L}$ such that $M_\beta \leftrightarrow {\tilde L}$.
    We show that the specified conditions on $[\gamma]$ must hold.

    If $i>j+1$, then $M_{\beta^{\rho}} \leftrightarrow {\bar
    L}$, where, by Lemma~\ref{le:groupoflm}, ${\bar L}$ is a
    solution to $\Ec_{i-1,j}(L)$. Then by the previous theorem the
    condition $[\gamma] = [\omega]^{\rho^{p-j}}$ is
    satisfied.

    If $i=j+1$ and $\Upsilon=1$, then let $[\alpha]\in J_1$ with
    $e([\alpha])\neq 0$. (See Lemma~\ref{le:smlgthnontriv}.
    Observe that $K(\root{p}\of{\alpha})/F$ is cyclic of degree
    $p^2$.) Since $\Ec'_{j+1,j}(L)$ is solvable there exists
    $M_\beta\supset M_\gamma,l(M_\beta)=j+1$ and $e([\beta])\neq 0$.
    Set $\beta'=\beta^{e([\alpha])}/\alpha^{e([\beta])}$. Then
    $e([\beta'])=0$. Since $\alpha\in J_1$ we have
    $M_{(\beta')^{\rho}}=M_{\beta^{\rho}}=M_\gamma$.
    By Proposition~\ref{pr:Kummer} and Lemma~\ref{le:groupoflm},
    $M_{\beta'}\leftrightarrow\bar{L}$, where $\bar{L}$ is a
    solution to $\Ec_{j+1,j}(L)$. By the previous theorem, the
    condition $[\gamma] = [\omega]^{\rho^{p-j}}$ is
    satisfied.

    Now consider the case $i=j+1$ and $\Upsilon=0$. By choosing
    another generator $[\beta]$ of $M_{\beta}$ if necessary, we may
    assume that $[\beta]^{\rho}=[\gamma]$. Consider $\beta'=\beta
    a^{-e([\beta])/p}$. Then $e([\beta'])=0$.  Now
    $[a^{-e([\beta])/p}]^{\rho}=[\xi_p^{-e([\beta])}]$, so
    $[\beta']^{\rho}=[\xi_p^{-e([\beta])}\gamma]$. Therefore if
    $l(M_{\xi_p^{-e([\beta])}\gamma})$ is at least $1$
    there exists a solution to an embedding problem corresponding
    to $M_{\beta'}$ and $M_{\xi_p^{-e([\beta])}\gamma}$. By the
    previous theorem, the condition $[\xi_p]^{-e([\beta])}
    [\gamma] = [\omega]^{\rho^{p-j}}$ is satisfied.
    If $[\xi_p^{-e([\beta])}\gamma]=[1]$ then we can set
    $\omega=1$ and again the condition
    $[\xi_p]^{-e([\beta])}[\gamma]=[\omega]^{\rho^{p-j}}$
    is satisfied.

    In all cases, then, we have shown that if $\Ec_{i,j}'(L)$ is
    solvable, the corresponding condition on $[\gamma]$ holds.

    Now suppose that the condition of part~\ref{it:ecij'1} holds:
    $[\gamma]=[\omega]^{\rho^{p-j}}$ for some $\omega$.
    Here we include the case $i=p$.

    If $\Upsilon=1$ then let $[\alpha]\in J_1$ with
    $e([\alpha])\neq 0$. (See Lemma~\ref{le:smlgthnontriv}.)
    Consider $\beta = \alpha\omega^{\rho^{p-i}}$.  If $i<p$,
    $e([\omega]^{\rho^{p-i}})=0$ and hence $e([\beta])\neq
    0$. Since $[\alpha]\in J_1$, $l(M_\beta) =
    l(M_{\omega^{\rho^{p-i}}}) = l(M_\gamma)+i-j=i$.
    Moreover, since $[\beta]^{\rho^{i-j}}=[\gamma]$,
    $M_\beta\supset M_\gamma$ and we have shown that
    $\Ec_{i,j}'(L)$ is solvable. In the case when $i=p$ and
    $\Upsilon=1$ we observed above that $\Ec_{p,j}'(L)$
    is equivalent with $\Ec_{p,j}(L)$ and also the solvability
    conditions are the same. Hence by Theorem~\ref{th:ecij} we see that
    $\Ec_{p,j}'(L)$ is solvable.

    If $\Upsilon=0$ and $i>j+1$, then let $\alpha=\root{p}\of{a}$.
    Consider $\beta=\alpha\omega^{\rho^{p-i}}$.  Again if $i<p$
    then $e([\omega]^{\rho^{p-i}})=0$ and hence $e([\beta])\neq
    0$. Since $[\alpha]\in J_2$ (see Lemma~\ref{le:smlgthnontriv})
    and $i>j+1$, $l(M_\beta) = l(M_{\omega^{\rho^{p-i}}}) =
    l(M_\gamma)+i-j=i$. Moreover, since $[\beta]^{\rho^{i-j}}=
    [\gamma]$, $M_\beta\supset M_\gamma$ and we have shown that
    $\Ec_{i,j}'(L)$ is solvable. (Observe that we employed the
    condition $i<p$ only to ensure that $e([\beta])\neq 0$ in this
    case. If $i=p$ then $e([\beta])$ plays no role, and therefore
    we have covered this case in the construction above.)

    Now suppose that the condition of part~\ref{it:ecij'2} holds:
    $[\gamma]=[\xi_p]^e[\omega]^{\rho^{p-j}}$ for some $\omega\in
    K^\times$ and $e\not\equiv 0 \bmod p$.  Let $\beta= a^{e/p}
    \omega^{\rho^{p-j-1}}$. If $j+1<p$ then because
    $e([\omega]^{\rho^{p-j-1}})=0$ we have $e([\beta]) = e
    \not\equiv 0 \bmod p$. Moreover, $[\beta]^{\rho}=
    [\xi_p^e][\xi_p^{-e}][\gamma]=[\gamma]$, so $M_\beta\supset
    M_\gamma$ and we have shown that $\Ec_{j+1,j}'(L)$ is
    solvable. Finally observe that if $j+1=p$ we showed at the
    beginning of our proof that both embedding problems
    $\Ec'_{p,p-1}$ and $\Ec_{p,p-1}$ are the same, and that also
    the conditions in Theorem~\ref{th:ecij} and
    Theorem~\ref{th:ecij'} for the existence of a solution of this
    problem are equivalent. Hence the existence of a solution in
    this case follows from Theorem~\ref{th:ecij}.

    Next we shall derive an explicit form of any solution field
    $\tilde{L}$ of our embedding problem.

    Observe that for any $f\in F^{\times}$ and $J_{p-1}\supset
    M_{\beta}\supsetneqq M_{\gamma}$, we have
    $l(M_{f\beta})=l(M_\beta)$, $M_{f\beta}\supset M_{\gamma}$,
    and $e([f\beta])=e([\beta])$. Recall that in the case
    $\Upsilon=1,\alpha$ is any element in $K^\times$ with
    $K(\root{p}\of{\alpha})/F$ cyclic of degree $p^2$, and in the
    case $\Upsilon=0,\alpha$ is $\root{p}\of{\alpha}$. By
    Proposition~\ref{pr:Kummer}, then,
    \begin{equation*}
        \tilde{L}=L(\root{p}\of{f\beta},
        \root{p}\of{\beta^{\rho}},
        \dots,\root{p}\of{\beta^{\rho^{p-1}}})
    \end{equation*}
    is a solution to the appropriate embedding problem for
    $\beta=\alpha\omega^{\rho^{p-i}}$ in the case $i>j+1-\Upsilon$
    or $j=p-1$ and $\beta=a^{e/p}\omega^{\rho^{p-i}}$ in the case
    $i=j+1$, $\Upsilon=0$, as above.

    To show that every solution field $\tilde{L}$ takes this form,
    suppose that $M_{\beta}\leftrightarrow\tilde{L}$ is a solution
    to $\Ec_{i,j}'(L)$. Hence $M_{\beta}\supset M_{\gamma}$. We
    consider first the case of part~\ref{it:ecij'1} when
    $\Upsilon=1$. If $i<p$ then by Lemma~\ref{le:groupoflm},
    $e([\beta])\neq 0$; in this case we let $c\in\mathbb{Z}$ be
    such that $e([\beta^c])=e([\alpha])$ and set $\beta'=\beta^c$
    so that $e([\beta'/\alpha])=0$. If $i=p$ then let
    $\beta'=\beta$. Because $i>j\ge 1$ and $[\alpha]\in J_1$,
    $l(M_{\beta'/\alpha})=l(M_{\beta'})=i$. Observe that
    $[\beta'/\alpha]^{\rho}=[\beta^c]^{\rho}$, so
    $M_{\beta'/\alpha}\supset M_{(\beta'/\alpha)^{\rho}} =
    M_{\beta^{\rho}}\supset M_\gamma$, because $M_\gamma$ is
    properly contained in $M_{\beta}$ and $M_{\beta^{\rho}}$ is
    the maximal proper $\Fp[G]$-submodule of $M_{\beta}$. Hence
    $M_{\beta'/\alpha}\leftrightarrow\bar{L}$, for $\bar{L}$ a
    solution to $\Ec_{i,j}(L)$.

    By Kummer theory and Theorem~\ref{th:ecij},
    \begin{equation*}
        M_{\beta'/\alpha}=\mo[f\omega^{\rho^{p-i}}],
        [\omega^{\rho^{p-i+1}}],\dots,
        [\omega^{\rho^{p-1}}]\mc
    \end{equation*}
    for some $f\in F^{\times}$ and $\omega\in K^{\times}$. Observe
    that hence $M_{\beta'/\alpha}=M_{f\omega^{\rho^{p-i}}}$.
    Because $[\beta'/\alpha]$ and $[f\omega^{\rho^{p-i}}]$ are
    both $\Fp[G]$-module generators of the same module of length
    $i$,
    \begin{equation*}
        \beta'/\alpha=(f\omega^{\rho^{p-i}})^{c_0 + c_1
        \rho+\dots+c_{i-1}\rho^{i-1}}
    \end{equation*}
    for some $c_k\in\Fp$. Let $f'=f^{c_0}$ and $\omega'=
    \omega^{c_0+c_1\rho+ \dots+c_{i-1}\rho^{i-1}}$. Then
    $\beta'/\alpha=f'(\omega')^{\rho^{p-i}}$, or
    $\beta'=f'\alpha(\omega')^{\rho^{p-i}}$. Since
    $M_{\beta'}=M_{\beta}\leftrightarrow\tilde{L}$, $\tilde{L}$
    takes the form
    \begin{equation*}
        \tilde{L}=L(\root{p}\of{f'\alpha(\omega')^{\rho^{p-i}}},
        \root{p}\of{(\omega')^{\rho^{p-i+1}}},
        \dots, \root{p}\of{(\omega')^{\rho^{p-1}}})
    \end{equation*}
    by Proposition~\ref{pr:Kummer}.

    The case when $\Upsilon=0$ can be treated as above with slight
    modifications. First in this case instead of $[\alpha]\in J_1$
    we take $[\root{p}\of{a}]$. We use our hypothesis $i>j+1$ to
    make sure as above that
    $l(M_{\beta'/\root{p}\of{a}})=l(M_{\beta'})$. Next observe
    that
    \begin{equation*}
        M_{\beta'/\root{p}\of{a}}\supset
        M_{(\beta'/\root{p}\of{a})^{\rho^2}}
        =M_{\beta^{\rho^{2}}}\supset M_\gamma
    \end{equation*}
    as the $\Fp[G]$-submodules of $M_\beta$ are linearly ordered
    by inclusion and $l(M_\beta)-l(M_\gamma)\geq 2$. The rest of
    the argument for case (1) when $\Upsilon=0$ faithfully follows
    the argument for case (1) when $\Upsilon=1$ as above.

    In order to show that every solution field $\tilde L$ in
    part~\ref{it:ecij'2} takes the specified form, observe that
    $A_j/A_{j+1}$ is in the center of $B_{j+1,e}$. Therefore,
    since we have one solution of the embedding problem
    $\Ec'_{j+1,j}$ of the form $L'= L(\root{p}\of{a^{e/p}
    \omega^{\rho^{p-j-1}}})$, by the well-known theorem on
    solutions of central embedding problems (see
    \cite[Lemma~A.1.1]{JLY}), any other solution $\tilde{L}$ of
    the embedding problem $\Ec'_{j+1,j}$ takes the form
    \begin{equation*}
        \tilde{L} = L(\root{p}\of{fa^{e/p}
        \omega^{\rho^{p-j-1}}}),\ \
        f\in F^{\times},
    \end{equation*}
    as required.
\end{proof}

\begin{remark*}
Lemma~\ref{le:surjections} implies that among our embedding
problems, only Galois extensions $L/F$ with $\Gal(L/F)\cong
B_{j,0}$ may be solved. This result agrees with our
Theorems~\ref{th:ecij} and~\ref{th:ecij'}, as follows. Suppose
that $L$ is the Galois closure of $K(\root{p}\of{\gamma})$,
$\gamma\in K^{\times}$, and $l(M_\gamma)=j$. From our solvability
conditions we see that if $L$ can be embedded in some extension
$\tilde{L}$ such that $\Gal(\tilde{L}/F)\cong B_{i,e}$ and $p\geq
i>j\geq 1,e\in\Fp$, then necessarily $e([\gamma])=0$.
\end{remark*}

\section{Arbitrary Fields}\label{se:arbitrary}

\subsection{Characteristic not $p$}\

We now suppose that $K_0/F_0$ is a cyclic extension of degree $p$
of fields of characteristic not $p$.  Set $F=F_0(\xi_p)$,
$K=K_0(\xi_p)$, and $s=[F:F_0]$.  Let $\epsilon$ denote a
generator of $\Gal(F/F_0)$ and $\sigma$ a generator of
$G=\Gal(K_0/F_0)$.  Since $p$ and $s$ are relatively prime,
$\Gal(K/F_0)\cong \Gal(F/F_0) \times \Gal(K_0/F_0)$. Therefore we
may naturally extend $\epsilon$ and $\sigma$ to $K$, and they
commute in $\Gal(K/F_0)$.  Using the extension of $\sigma$ to $K$,
we write $G$ for $\Gal(K/F)$ as well.

Let $t\in \Z$ such that $\epsilon(\xi_p)=\xi_p^t$. Then $t$ is
relatively prime to $p$.  Let $J^\epsilon$ be the $t$-eigenspace
of $J=K^\times/K^{\times p}$ under the action of $\epsilon$.
Observe that since $\epsilon$ and $\sigma$ commute, $J^\epsilon$
is an $\Fp[G]$-subspace of $J$. By \cite[\S 5, Prop.]{Wat}, we
have a Kummer correspondence over $K_0$ of subspaces $M^\epsilon$
of the $\Fp$-vector space $J^\epsilon$ and abelian exponent $p$
extensions $L_0$ of $K_0$:
\begin{multline*}
    M^\epsilon = ((KL_0)^{\times p} \cap K^\times)/K^{\times p}
    \quad \leftrightarrow \\ L_0=\text{maximal $p$-extension
    of $K_0$ in\ } L_{M^\epsilon}=K(\root{p}\of{\gamma}:[\gamma]
    \in M_\epsilon).
\end{multline*}
As Waterhouse shows, for $M^\epsilon\subset J^\epsilon$,
$\epsilon\in \Gal(K/K_0)$ has a unique lift $\tilde\epsilon$ to
$\Gal(L_{M^\epsilon}/K_0)$ of order $s$, and $L_0$ is the fixed
field of $\tilde\epsilon$.

We first prove a lemma on the decomposition of $J$:

\begin{lemma}\label{le:jdecomp}
    $J=J^\epsilon\oplus J^\nu$, where $J^\nu$ is an
    $\Fp[G]$-submodule of $J$, and $e$ is trivial on
    $J^\nu \cap J_{p-1}$.
\end{lemma}

\begin{proof}
We adapt an approach to descent from \cite[page 258]{Sa}. Let
$z\in \Z$ satisfy $zst^{s-1}\equiv 1 \bmod p$, and set
\begin{equation*}
    T = z\cdot \sum_{i=1}^st^{s-i}\epsilon^{i-1} \in
    \Z[\Gal(K/F_0)].
\end{equation*}
We calculate that $(t-\epsilon)T\equiv 0 \mod p$, and hence the
image of $T$ on $J$ lies in $J^\epsilon$.  Moreover, on
$J^\epsilon$, $\epsilon$ acts by multiplication by $t$, and hence
$T$ acts as the identity on $J^\epsilon$.  Finally, since
$\epsilon$ and $\sigma$ commute, $T$ and $I-T$ commute with
$\sigma$. Hence $J$ decomposes into a direct sum $J^\epsilon\oplus
J^\nu$, with associated projections $T$ and $I-T$.

Let $a\in F^\times$ satisfy $K=F(\root{p}\of{a})$, and consider
$[a]_F\in F^\times/F^{\times p}$. By \cite[\S 5, Prop.]{Wat},
$\epsilon([a]_F)=[a]_F^t$.  Suppose $\gamma\in K^\times$ satisfies
$[\gamma]\in J_{p-1}$.  Then, since $\epsilon$ and $\sigma$
commute,
\begin{equation*}
    [N(\epsilon(\gamma))]_F = [\epsilon(N(\gamma))]_F =
    \epsilon([N(\gamma)]_F) = [N(\gamma)]_F^t.
\end{equation*}
Hence $e(\epsilon([\gamma]))=t\cdot e([\gamma])$, and we then
calculate that $e([T\gamma])=e([\gamma])$. Therefore
$e((I-T)[\gamma])=0$.
\end{proof}

Now we establish that the Galois structure of $L_{M^\epsilon}/F$
is equivalent to that of $L_0/F_0$.

\begin{proposition}\label{pr:Kummer2}
    Under the Kummer correspondence above, $L_0$ is Galois over
    $F_0$ if and only if $M^\epsilon$ is an $\Fp[G]$-submodule of
    $J^\epsilon$.  In this case the base extension $F_0\to F$ induces
    a natural isomorphism of $G$-extensions $\Gal(L_0/F_0)\cong
    \Gal(L/F)$.
\end{proposition}

\begin{proof}
    If $L_0/F_0$ is Galois, then $L_{M^\epsilon}=L_0K/F$ is Galois
    as well, and by Proposition~\ref{pr:Kummer} (\ref{item1}),
    $L^\epsilon$ is an $\Fp[G]$-submodule of $J$.

    Going the other way, suppose that $M^\epsilon$ is an
    $\Fp[G]$-submodule of $J^\epsilon$.  By the correspondence,
    $L_{M^\epsilon}/K_0$ is Galois. Then $M^\epsilon$ is also an
    $\Fp[\Gal(K/F_0)]$-submodule of $J^\epsilon$ and therefore
    $L_{M^\epsilon}/F_0$ is Galois.

    Now since $K_0/F_0$ is Galois, every automorphism of
    $L_{M^\epsilon}$ sends $K_0$ to $K_0$.  Moreover,
    since $L_0$ is the unique maximal $p$-extension of $K_0$ in
    $L_{M^\epsilon}$, every automorphism of $L_{M^\epsilon}$ sends
    $L_0$ to $L_0$.  Therefore $L_0/F_0$ is Galois.

    Finally we show that base field extension $F_0\to F$ induces
    a natural isomorphism of $G$-extensions $\Gal(L_0/F_0)\to
    \Gal(L/F)$.  Now $F/F_0$ and $L_0/F_0$ are of relatively prime
    degrees, and hence $\Gal(L_0F/F_0)\cong \Gal(F/F_0)\times
    \Gal(L_0/F_0)$. Moreover we deduce that we have the natural
    isomorphism $G=\Gal(K_0/F_0)\cong \Gal(K/F)$, and that the natural isomorphism
    $\Gal(L/F)\cong\Gal(L_0/F_0)$ is a $G$-extension isomorphism.
\end{proof}

Now we make the connection between embedding problems over $F_0$
and embedding problems over $F$.  For $p\ge i>j\ge 1$, denote by
$\Ec_{i,j}(L_0)$ and $\Ec'_{i,j}(L_0)$ the embedding problems
\begin{equation*}
    \Ec_{i,j}(L_0)\colon \quad 1\to A_j/A_i \to B_{i,0} \to
    (A/A_j)\rtimes G=\Gal(L_0/F_0)\to 1
\end{equation*}
and, for any $e\neq 0$,
\begin{equation*}
    \Ec_{i,j}'(L_0)\colon \quad 1\to A_j/A_i \to B_{i,e}\to
    (A/A_j)\rtimes G=\Gal(L_0/F_0)\to 1.
\end{equation*}

\noindent {\bf Warning.} \

In order to avoid possible confusion,
let us recall that by $[\delta]^\epsilon$ we mean the projection
of $[\delta]$ into the summand $J^\epsilon$ of $J$. Similarly
$[\delta]^\nu$ means the projection of $[\delta]$ into
the summand $J^\nu$ of $J$.

\begin{proposition}\label{pr:equiv}\
    \begin{enumerate}
        \item $\Ec_{i,j}(L_0)$ is solvable if and only if
        $\Ec_{i,j}(L)$ is solvable.
        \item $\Ec'_{i,j}(L_0)$ is solvable if and only if
        $\Ec'_{i,j}(L)$ is solvable.
    \end{enumerate}
\end{proposition}

\begin{proof}
    Let $\tilde L_0$ be a solution to $\Ec_{i,j}(L_0)$.  Then by
    Proposition~\ref{pr:Kummer2}, $\tilde L:=\tilde L_0F$ is a
    solution to $\Ec_{i,j}(L)$.

    Going the other way, let $\tilde L$ be a solution to
    $\Ec_{i,j}(L)$.  Kummer theory gives correspondences
    $M^\epsilon \leftrightarrow L_0$ over $K_0$, as well as
    $M^\epsilon\leftrightarrow L$ and $\tilde M\leftrightarrow
    \tilde L$ over $K$.  By Proposition~\ref{pr:Kummer}
    (\ref{item2}), $\tilde M = M_\delta$ and $M^\epsilon=M_\gamma$
    for some $\delta, \gamma\in K^\times$, with $[\gamma]\in
    J^\epsilon$. By Lemma~\ref{le:jdecomp}, we may write $[\delta]
    = [\delta]^\epsilon+[\delta]^\nu\in J^\epsilon \oplus J^\nu$,
    with $e([\delta]^\epsilon) = e([\delta])$ if $[\delta]\in J_{p-1}$.
    Moreover, since $[\delta]^{\rho^{i-j}}=[\gamma]$ and $J^\nu$ is a
    $\Fp[G]$-submodule, $([\delta]^{\epsilon})^ {\rho^{i-j}} =
    [\gamma]$.  Let $\tilde M^\epsilon = M_{\delta^\epsilon}$.
    Then $M^\epsilon\subset \tilde M^\epsilon$.  By the Kummer
    correspondence over $K_0$, there exists a field $\tilde L_0$
    such that $\tilde M^\epsilon \leftrightarrow \tilde L_0$ and
    $L_0\subset \tilde L_0$.  By Lemma~\ref{le:groupoflm} and
    Proposition~\ref{pr:Kummer2}, $\tilde L_0$ is a solution to
    $\Ec_{i,j}(L_0)$.

    The case of $\Ec'_{i,j}$ follows analogously.
\end{proof}

\subsection{Embedding Problem Conditions, Arbitrary Fields}\

To state the general result, we alter our notation to let $F$ take
the place of $F_0$.  If $\chr F\neq p$, then $J$ now denotes
$K(\xi_p)^{\times}/K(\xi_p)^{\times p}$.

\begin{theorem}\label{th:arb}
    Let $F$ be an arbitrary field.
    \begin{enumerate}
        \item If $\chr F=p$, then $\Ec_{i,j}(L)$ and
        $\Ec'_{i,j}(L)$ are solvable.
        \item If $\chr F\neq p$, let $\gamma\in F(\xi_p)^\times$
        satisfy $K(\xi_p)=F(\xi_p)(\root{p}\of{\gamma})$, and set
        $\Upsilon=1$ if $\xi_p\in N_{K(\xi_p)/F(\xi_p)}
        (K(\xi_p)^\times)$ and $\Upsilon=0$ otherwise.  Then
        \begin{enumerate}
        \item $\Ec_{i,j}(L)$ is solvable if and only if
        $[\gamma]=[\omega]^{\rho^{p-j}}$ in $J$ for some
        $\omega\in K(\xi_p)^\times$.
        \item $\Ec'_{i,j}(L)$, $i>j+1-\Upsilon$ or $j=p-1$, is
        solvable if and only if $[\gamma]=[\omega]^{\rho^{p-j}}$
        in $J$ for some $\omega\in K(\xi_p)^\times$.
        \item $\Ec'_{j+1,j}(L)$, $\Upsilon=0$, is solvable if and
        only if $[\gamma]=[\xi_p]^e[\omega]^{\rho^{p-j}}$ in $J$
        for some $\omega\in K(\xi_p)^\times$ and $e\not\equiv 0
        \bmod p$.
        \end{enumerate}
    \end{enumerate}
\end{theorem}

\begin{proof}
    If $\chr F=p$, then by Witt's theorem all central non-split embedding
    problems with kernel $\Fp$ are solvable. (See \cite[Appendix~A]{JLY}.)
    Since $B_{i,e}$ and $B_{j,0}$ have the same minimal
    number of generators for all $1\le i,j$ and $e$, Witt's
    theorem gives that $\Ec_{i,j}(L)$ and $\Ec'_{i,j}(L)$ are
    solvable. Indeed one can successively solve a chain of
    suitable central non-split embedding problems with kernel $\Fp$ leading
    to solutions of $\Ec_{i,j}(L)$ and $\Ec'_{i,j}(L)$.

    If $\chr F\neq p$, then the statements follow from
    Theorems~\ref{th:ecij} and \ref{th:ecij'} using
    Proposition~\ref{pr:equiv}.
\end{proof}

\section{Proof of Main Theorem}\label{se:prfmain}

\begin{proof}
    Observe that $\Ec_i=\Ec_{i,1}$.  The equivalence of
    (\ref{basicitem1}) and (\ref{basicitem2}) follows from
    Theorem~\ref{th:arb}.

    Now assume that $\chr F\neq p$.  By
    Proposition~\ref{pr:equiv}, $\Ec_{i,1}(L)$ is solvable if and
    only if $\Ec_{i,1}(L(\xi_p))$ is solvable. The condition on
    $b$ implies that $M_b\leftrightarrow L(\xi_p)/K(\xi_p)$ under
    the Kummer correspondence.  To show (\ref{basicitem3}), by
    Theorem~\ref{th:ecij} we need only show that there exists
    $\alpha\in K(\xi_p)^\times$ with $[b]=[\alpha]^ {\rho^{p-1}}$
    if and only if there exists $\omega$ with
    $N_{K(\xi_p)/F(\xi_p)}(\omega)=b$.  Let $N$ denote
    $N_{K(\xi_p)/F(\xi_p)}$.

    By equation~\eqref{eq:equiv} of Section~\ref{se:index},
    \begin{equation*}
        [\omega]^{\rho^{p-1}} = [\omega]^{1+ \sigma+ \dots+
        \sigma^{p-1}} = [\omega^{1+ \sigma+ \dots+ \sigma^{p-1}}]
        = [N(\omega)]=[b],
    \end{equation*}
    so if $\omega$ exists satisfying $N(\omega)=b$, then $\alpha =
    \omega$ satisfies $[b]=[\alpha]^ {\rho^{p-1}}$.

    Going the other way, suppose that $[b] = [\alpha]^
    {\rho^{p-1}}$ for some $\alpha$.  By equation~\eqref{eq:equiv}
    again, $[b]=[N(\alpha)]$.  Then $N(\alpha)=k^pb$ for some
    $k\in K(\xi_p)^\times$.  Hence $N(\alpha)/b\in F(\xi_p)^\times
    \cap K(\xi_p)^{\times p}$.  Let $a\in K(\xi_p)^\times$ satisfy
    $K(\xi_p)=F(\xi_p)(\root{p}\of{a})$. Then by Kummer theory
    $k^p=a^s f^p$ for some $s\in \Z$. Choosing
    $\omega=\alpha/(a^{s/p}f)$ we obtain $N(\omega)=b$.

    Finally, the explicit solution of solution fields in the
    case $\xi_p\in F$ follows directly from Theorem~\ref{th:ecij}.
\end{proof}

\section*{Acknowledgements}

It is clear from the content of our paper that we are strongly
influenced by the beautiful paper of Waterhouse cited in the
references, and we gratefully acknowledge this crucial influence
on our work. We are also grateful to Adrian Wadsworth for some
helpful discussions concerning this work.  Finally, we thank an
anonymous referee for suggesting that our results would likely
extend from the case of fields with a primitive $p$th root
of unity to the arbitrary field case, using Albert's descent technique.

\end{document}